\documentclass{article}
\usepackage{amsfonts}
\usepackage{amsmath}
\usepackage{amsmath,amssymb,latexsym,color}
\usepackage[mathscr]{eucal}
\newtheorem{thm}{Theorem}[section]
\newtheorem{prop}[thm]{Proposition}
\newtheorem{lemma}[thm]{Lemma}
\newtheorem{cor}[thm]{Corollary}

\newtheorem{example}{Example}[section]
\newtheorem{defin}[thm]{Definition}

\newtheorem{remark}{Remark}[section]

\newcommand{\qed}{\hfill\Box\medskip}
\usepackage{CJK}
\begin{document}
\begin{CJK*}{GBK}{song}

\renewcommand{\baselinestretch}{1.3}
\title{\bf The Terwilliger algebra of Odd graphs}

\author{
Qian Kong\quad Benjian Lv\quad Kaishun Wang\footnote{ Corresponding
author.\newline \quad\quad E-mail addresses:
kongqian@mail.bnu.edu.cn (Qian Kong);
 benjian@mail.bnu.edu.cn (Benjian Lv); wangks@bnu.edu.cn (Kaishun Wang)} \\
{\footnotesize \em  Sch. Math. Sci. {\rm \&} Lab. Math. Com. Sys.,
Beijing Normal University, Beijing, 100875,  China} }
\date{}
\maketitle

\begin{abstract}
In [The Terwilliger algebra of the Johnson schemes, Discrete
Mathematics 307 (2007) 1621--1635], Levstein and Maldonado computed
the Terwilliger algebra of the Johnson scheme $J(n,m)$ when $3m\leq
n$. The distance-$m$ graph of $J(2m+1,m)$ is the Odd graph
$O_{m+1}$. In this paper, we determine the Terwilliger algebra of
$O_{m+1}$ and give its basis.

\medskip
\noindent {\em AMS classification:} 05E30

\noindent {\em Key words:} Terwilliger algebra; Odd graph

\end{abstract}

\section{Introduction}

Suppose $\Gamma=(X,R)$ denotes a simple connected graph with
diameter $D$. For each $i\in \{0,1,\ldots,D\}$, let
$\Gamma_i(x)=\{y\in X\mid
\partial(x,y)=i\}$, where $\partial(x,y)$ is the distance between $x$ and $y$. Define $E_i^*=E_i^*(x)$ to be the diagonal matrix in
Mat$_X(\mathbb{C})$ with $yy$-entry
$$
(E_i^*)_{yy}=\left\{ \begin{array}{ll}
1,& \textrm{if $y\in \Gamma_i(x)$},\\
0,& \textrm{otherwise}.\\
\end{array} \right.
$$
The \emph{Terwilliger algebra} $\mathcal {T}(x)$ of $\Gamma$ with
respect to a given vertex $x$ is the subalgebra of
Mat$_X(\mathbb{C})$ generated  by the adjacency matrix $A$ of
$\Gamma$ and $E_0^*,E_1^*,\ldots,E_D^*$.

Terwilliger \cite{ter} initiated the study of the Terwilliger
algebra of association schemes, which has been used to study
(almost) bipartite $P$- and $Q$-polynomial association schemes
\cite{bpq, abpq}, $2$-homogeneous bipartite distance-regular graphs
\cite{curtin}, Hypercubes \cite{go}, Hamming graphs \cite{hdq},
Johnson graphs \cite{johnson}, incidence graphs of Johnson geometry
\cite{kong} and so on.

Let $\Omega$ be a set of cardinality $2m+1$ and let ${\Omega \choose
i}$ denote the set of all $i$-subsets of $\Omega$. The Odd graph
$O_{m+1}$ is the graph whose vertex set is the set $X={\Omega
\choose m}$, where two vertices are adjacent if they are disjoint.
Levstein and Maldonado \cite{johnson} determined the Terwilliger
algebra of the Johnson graph $J(n,m)$ when $3m\leq n$. Observe
$O_{m+1}$ is the distance-$m$ graph of the Johnson graph
$J(2m+1,m)$, and they have the same Terwilliger algebra. In this
paper we shall determine the Terwilliger algebra of $O_{m+1}$
(Theorem \ref{17}), give one of its bases (Proposition \ref{18}) and
compute its dimension (Corollary \ref{21}).

\section{Intersection matrix}

In this section we first introduce the intersection matrix, then
discover the relationship between the adjacency matrix of the Odd
graph $O_{m+1}$ and the intersection matrices.

Since $O_{m+1}$ is distance-transitive with diameter $m$ (cf. \cite
{bcn}), the isomorphism class of $\mathcal {T}(x)$ is independent of
the choice of $x$, denoted by $\mathcal {T}:=\mathcal {T}(x)$.

Let $V$ be a set of cardinality $v$. Let $H_{i,j}^l(v)$ be a binary
matrix with rows indexed by ${V \choose i}$ and columns indexed by
$V \choose j$, whose
 $yz$-entry is defined by
$$
(H_{i,j}^l(v))_{yz}=\left\{ \begin{array}{ll}
1,& \textrm{if $|y\cap z|=l$},\\
0,& \textrm{otherwise}.\\
\end{array} \right.
$$
This matrix is a class of \emph{intersection matrices}. Observe that
$H_{i,j}^l(v)\neq 0$ if and only if $\max(0,i+j-v)\leq l\leq
\min(i,j)$. We adopt the convention that $H_{i,j}^l(v)=0$ for any
integer $l$ such that $l<0$ or $l>\min(i,j)$. From \cite
[Proposition 4] {mm}, we have
\begin{eqnarray}\label{7}
H_{i,j}^l(v)H_{j,k}^s(v)=\sum_{g=0}^{\min(i,k)}\sum_{h=0}^g{g
\choose h}{i-g \choose l-h}{k-g \choose s-h}{v+g-i-k \choose
j+h-l-s}H_{i,k}^g(v).
\end{eqnarray}
In particular,
\begin{eqnarray}\label{8}
H_{i,j}^l(v)H^0_{j,k}(v)=\sum_{s=\max(0,i+j+k-l-v)}^{\min(i-l,k)}{i-s
\choose l}{v+s-i-k \choose j-l}H_{i,k}^s(v).
\end{eqnarray}

\begin{lemma}\label{2}
Let $\Gamma$ be the Odd graph $O_{m+1}$ with the adjacency matrix
$A$, and let $A_{i,j}$ be the submatrix of $A$ with rows indexed by
$\Gamma_i(x)$ and columns indexed by $\Gamma_j(x)$. Then
\begin{eqnarray}
\label{3}{}&&A_{i,j}=0 \quad (0\leq i\leq j \leq m, \ i\neq j-1 \
\textrm{or} \ i=j\neq m),\\
\label{4}{}&&A_{2i,2i+1}=H^0_{m-i,i}(m)\otimes H^0_{i,m-i}(m+1) \quad (0\leq i\leq \lceil\frac{m}{2}\rceil-1),\\
\label{5}{}&&A_{2i+1,2i+2}=H^0_{i,m-i-1}(m)\otimes
H^0_{m-i,i+1}(m+1) \quad (0\leq i\leq
\lfloor\frac{m}{2}\rfloor-1),\\
\label{6}{}&&A_{m,m}=H^0_{\lfloor\frac{m}{2}\rfloor,\lfloor\frac{m}{2}\rfloor}(m)\otimes
H^0_{\lceil\frac{m}{2}\rceil,\lceil\frac{m}{2}\rceil}(m+1),
\end{eqnarray}
where ``$\otimes$'' denotes the Kronecker product of matrices.
\end{lemma}
\textbf{Proof.} Since $O_{m+1}$ is almost bipartite, (\ref{3}) is
directed.

Pick $y\in\Gamma_{2i}(x)$, $z\in\Gamma_{2i+1}(x)$. Note that
$\partial(x,y)=2i$ if and only if $|x\cap y|=m-i$;
$\partial(x,z)=2i+1$ if and only if $|x\cap z|=i$. Then $|x\cap
y|=m-i$ and $|x\cap z|=i$. Suppose
$y=\alpha_{m-i}\beta_i:=\alpha_{m-i}\cup\beta_i$,
$z=\alpha_i'\beta_{m-i}'$, where $\alpha_{m-i}\in{x \choose m-i}$
and $\beta_i\in{\Omega\setminus x \choose i}$, while $\alpha_i'\in{x
\choose i}$ and $\beta_{m-i}'\in{\Omega\setminus x \choose m-i}$.
Then
$$
(A_{2i,2i+1})_{yz}=(H^0_{m-i,i}(m)\otimes
H^0_{i,m-i}(m+1))_{yz}=\left\{ \begin{array}{ll}
1,& \textrm{if $\alpha_{m-i}\cap\alpha_i'=\emptyset$ and $\beta_i\cap\beta_{m-i}'=\emptyset$},\\
0,& \textrm{otherwise},\\
\end{array} \right.
$$
which leads to (\ref{4}).

Similarly, (\ref{5}) (\ref{6}) hold. $\qed$

\section{The Terwilliger algebra}

In this section we fix $x\in {\Omega \choose m}$, then consider the
Terwilliger algebra $\mathcal {T}=\mathcal {T}(x)$ of $O_{m+1}$.

For $0\leq i,j\leq m$, any matrix $M$ indexed by elements in
$\Gamma_i(x)\times\Gamma_j(x)$ can be embedded into
Mat$_X(\mathbb{C})$ by
$$
L(M)_{\Gamma_p(x)\times\Gamma_q(x)}=\left\{ \begin{array}{ll}
M,& \textrm{if $p=i$ and $q=j$},\\
0,& \textrm{otherwise}.\\
\end{array} \right.
$$
Write $G_{i,j}(v)=\{g\mid \max(0,i+j-v)\leq g\leq \min(i,j)\}$. Let
\begin{eqnarray}\label{9}
\mathcal {M}=\bigoplus\limits_{p,q=0}^{m}L(\mathcal {M}_{p,q}),
\end{eqnarray}
where $L(\mathcal {M}_{p,q})=\{L(M)\mid M\in \mathcal {M}_{p,q}\}$,
and
\begin{eqnarray}\label{1}
\mathcal {M}_{2i,2j}=\textrm{Span}\{H_{m-i,m-j}^l(m)\otimes
H_{i,j}^s(m+1),\ l\in G_{m-i,m-j}(m), \ s\in G_{i,j}(m+1)\},
\end{eqnarray}
\begin{eqnarray}\label{19}
\mathcal {M}_{2i,2j+1}=\textrm{Span}\{H_{m-i,j}^l(m)\otimes
H_{i,m-j}^s(m+1),\ l\in G_{m-i,j}(m), \ s\in G_{i,m-j}(m+1)\},
\end{eqnarray}
\begin{eqnarray}\label{22}
\mathcal {M}_{2i+1,2j}=\textrm{Span}\{H_{i,m-j}^l(m)\otimes
H_{m-i,j}^s(m+1),\ l\in G_{i,m-j}(m), \ s\in G_{m-i,j}(m+1)\},
\end{eqnarray}
\begin{eqnarray}\label{23}
\mathcal {M}_{2i+1,2j+1}=\textrm{Span}\{H_{i,j}^l(m)\otimes
H_{m-i,m-j}^s(m+1),\ l\in G_{i,j}(m), \ s\in G_{m-i,m-j}(m+1)\}.
\end{eqnarray}

Note that $\mathcal {M}$ is a vector space. By (\ref{7}) we have
$\mathcal {M}$ is an algebra. Next we shall prove $\mathcal
{T}=\mathcal {M}$.

\begin{lemma}\label{10}
The Terwilliger algebra $\mathcal {T}$ is a subalgebra of $\mathcal
{M}$.
\end{lemma}
\textbf{Proof.} By Lemma \ref{2} we have $ A\in \mathcal {M}$. Since
$$
E_{2i}^*=L(H_{m-i,m-i}^{m-i}(m)\otimes H_{i,i}^i(m+1))\in\mathcal
{M}, \quad 0\leq i\leq\lfloor\frac{m}{2}\rfloor,
$$
$$
E_{2i+1}^*=L(H_{i,i}^i(m)\otimes H_{m-i,m-i}^{m-i}(m+1))\in\mathcal
{M}, \quad 0\leq i\leq\lceil\frac{m}{2}\rceil-1,
$$
we have $\mathcal {T}\subseteq\mathcal {M}$. $\qed$

For $0\leq i,j\leq m$, let $\mathcal {T}_{i,j}=\{M_{i,j}\mid
M\in\mathcal {T}\}$, where $M_{i,j}$ is the submatrix of $M$ with
rows indexed by $\Gamma_i(x)$ and columns indexed by $\Gamma_j(x)$.
Since $\mathcal {T}$ is an algebra, each $\mathcal {T}_{i,j}$ is a
linear space. From $\mathcal {T}E_j^*\mathcal {T}\subseteq\mathcal
{T}$ we obtain $(\mathcal {T}E_j^*\mathcal
{T})_{i,k}\subseteq\mathcal {T}_{i,k}$, which implies that
\begin{eqnarray}\label{11}
\mathcal {T}_{i,j}\mathcal {T}_{j,k}\subseteq\mathcal {T}_{i,k}.
\end{eqnarray}
From $A$, $E_i^*\in\mathcal {T}$, we have
$AE_{i_2}^*AE_{i_3}^*\cdots AE_{i_{n-1}}^*A\in\mathcal {T}$, which
follows that
\begin{eqnarray}\label{12}
A_{i_1,i_2}A_{i_2,i_3}\cdots A_{i_{n-2},i_{n-1}}A_{i_{n-1},i_n}\in
\mathcal {T}_{i_1,i_n}.
\end{eqnarray}

\begin{lemma}\label{24}
For $i\leq j\leq \lceil\frac{m}{2}\rceil-1$ and $0\leq l\leq i$, we
have
$$H_{i,j}^l(m)\otimes H_{m-i,m-j}^{m-j}(m+1)\in\mathcal
{T}_{2i+1,2j+1}.
$$
\end{lemma}
\textbf{Proof.} We use induction on $l$.

By (\ref{12}), for $i<j$ we have $A_{2i+1,2i+2}A_{2i+2,2i+3}\cdots
A_{2j,2j+1}\in\mathcal {T}_{2i+1,2j+1}$, which yields that
\begin{eqnarray}\label{14}
H_{i,j}^i(m)\otimes H_{m-i,m-j}^{m-j}(m+1)\in\mathcal
{T}_{2i+1,2j+1}.
\end{eqnarray}
When $i=j$ we pick $I_{m \choose i}\otimes I_{m+1 \choose
m-i}\in\mathcal {T}_{2i+1,2j+1}$, which also satisfies (\ref{14}).

Assume that $H_{i,j}^g(m)\otimes H_{m-i,m-j}^{m-j}(m+1)\in\mathcal
{T}_{2i+1,2j+1}$ for $g\geq l$. By (\ref{11}) and (\ref{12}), for
$2j+1<m$ we obtain
$$
(H_{i,j}^l(m)\otimes
H_{m-i,m-j}^{m-j}(m+1))A_{2j+1,2j+2}A_{2j+2,2j+1}\in \mathcal
{T}_{2i+1,2j+1}\mathcal {T}_{2j+1,2j+1}\subseteq\mathcal
{T}_{2i+1,2j+1},
$$
and for $2j+1=m$  we have
$$
(H_{i,j}^l(m)\otimes
H_{m-i,m-j}^{m-j}(m+1))A_{2j+1,2j+1}^2\in\mathcal {T}_{2i+1,2j+1}.
$$
Then we get
$$
(a_1H_{i,j}^{l-1}(m)+a_2H_{i,j}^l(m)+a_3H_{i,j}^{l+1}(m))\otimes
H_{m-i,m-j}^{m-j}(m+1)\in\mathcal {T}_{2i+1,2j+1},
$$
where $a_1$, $a_2$ and $a_3$ are some positive integers. It follows
that $H_{i,j}^{l-1}(m)\otimes H_{m-i,m-j}^{m-j}(m+1)\in\mathcal
{T}_{2i+1,2j+1}$. Hence the conclusion is obtained by induction.
$\qed$

\begin{lemma}\label{25}
For $i+1\leq j\leq \lfloor\frac{m}{2}\rfloor$ and $0\leq l\leq i$,
we have
$$
H_{i,m-j}^l(m)\otimes H_{m-i,j}^{j-i-1}(m+1)\in\mathcal
{T}_{2i+1,2j}.
$$
\end{lemma}
\textbf{Proof.} By (\ref{12}) we have
$A_{2i+1,2i+2}A_{2i+2,2i+3}\cdots A_{2j-1,2j}\in\mathcal
{T}_{2i+1,2j}$ for $i+1\leq j$, i.e.,
$$
H_{i,m-j}^0(m)\otimes H_{m-i,j}^{j-i-1}(m+1)\in\mathcal
{T}_{2i+1,2j}.
$$

Assume that $H_{i,m-j}^g(m)\otimes H_{m-i,j}^{j-i-1}(m+1)\in\mathcal
{T}_{2i+1,2j}$ for $g\leq l$. Then by (\ref{11}) and (\ref{12}), we
obtain
$$
(H_{i,m-j}^l(m)\otimes
H_{m-i,j}^{j-i-1}(m+1))A_{2j,2j-1}A_{2j-1,2j}\in\mathcal
{T}_{2i+1,2j},
$$
which gives
$$
(b_1H_{i,m-j}^{l-1}(m)+b_2H_{i,m-j}^l(m)+b_3H_{i,m-j}^{l+1}(m))\otimes
H_{m-i,j}^{j-i-1}(m+1)\in\mathcal {T}_{2i+1,2j},
$$
where $b_1$, $b_2$ and $b_3$ are some positive integers. Thus
$H_{i,m-j}^{l+1}(m)\otimes H_{m-i,j}^{j-i-1}(m+1)\in\mathcal
{T}_{2i+1,2j}$ and the conclusion is valid by induction. $\qed$

\begin{lemma}\label{20}
The algebra $\mathcal {M}$ is a subalgebra of $\mathcal {T}$.
\end{lemma}
\textbf{Proof.} In order to prove this result, we only need to show
that $\mathcal {M}_{p,q}\subseteq\mathcal {T}_{p,q}$ for $0\leq
p,q\leq m$. Write $\mathcal {M}_{p,q}^\textrm{t}=\{M^\textrm{t}\mid
M\in\mathcal {M}_{p,q}\}$ and $\mathcal
{T}_{p,q}^\textrm{t}=\{M^\textrm{t}\mid M\in\mathcal {T}_{p,q}\}$  .
Since $\mathcal {M}_{q,p}=\mathcal {M}_{p,q}^\textrm{t}$ and
$\mathcal {T}_{q,p}=\mathcal {T}_{p,q}^\textrm{t}$, it suffices to
prove $\mathcal {M}_{p,q}\subseteq\mathcal {T}_{p,q}$ for $p\leq q$.
We use induction on $p$.

\medskip
\textbf{Step 1.} Show $\mathcal {M}_{0,q}\subseteq\mathcal
{T}_{0,q}$ $(0\leq q\leq m)$.

According to (\ref{1}) (\ref{19}), we get
$$
\mathcal {M}_{0,2j}=\textrm{Span}\{H_{m,m-j}^{m-j}(m)\otimes
H_{0,j}^0(m+1)\} \quad (0\leq j\leq\lfloor\frac{m}{2}\rfloor),
$$
and
$$
\mathcal {M}_{0,2j+1}=\textrm{Span}\{H_{m,j}^j(m)\otimes
H_{0,m-j}^0(m+1)\} \quad (0\leq j\leq\lceil\frac{m}{2}\rceil-1).
$$
By Lemma \ref{2} and (\ref{8}), we have
$$
A_{0,1}A_{1,2}\cdots A_{2j-1,2j}=c_1H_{m,m-j}^{m-j}(m)\otimes
H_{0,j}^0(m+1)
$$
and
$$
A_{0,1}A_{1,2}\cdots A_{2j,2j+1}=c_2H_{m,j}^j(m)\otimes
H_{0,m-j}^0(m+1),
$$
where $c_1$, $c_2$ are some positive integers. Then by (\ref{12}) we
have $\mathcal {M}_{0,2j}\subseteq\mathcal {T}_{0,2j}$ and $\mathcal
{M}_{0,2j+1}\subseteq\mathcal {T}_{0,2j+1}$.

\medskip
\textbf{Step 2.} Assume that $\mathcal {M}_{p,q}\subseteq\mathcal
{T}_{p,q}$ for $p\leq 2i$. We will show that $\mathcal
{M}_{2i+1,q}\subseteq\mathcal {T}_{2i+1,q}$ and $\mathcal
{M}_{2i+2,q}\subseteq\mathcal {T}_{2i+2,q}$.

\medskip
\textbf{Step 2.1.} Show $\mathcal {M}_{2i+1,q}\subseteq\mathcal
{T}_{2i+1,q}$ $(2i+1\leq q\leq m)$.

\medskip
\textbf{Case 1.} $q=2j+1$ ($i\leq j\leq\lceil\frac{m}{2}\rceil-1$).

By inductive hypothesis we have
$$
H_{m-i,j}^l(m)\otimes H_{i,m-j}^s(m+1)\in\mathcal
{M}_{2i,2j+1}\subseteq\mathcal {T}_{2i,2j+1}, \quad l\in
G_{m-i,j}(m),\ s\in G_{i,m-j}(m+1).
$$
Since $A_{2i+1,2i}=A_{2i,2i+1}^\textrm{t}=H^0_{i,m-i}(m)\otimes
H^0_{m-i,i}(m+1)\in\mathcal {T}_{2i+1,2i}$, by (\ref{11}) we have
$$
(H^0_{i,m-i}(m)\otimes H^0_{m-i,i}(m+1))(H_{m-i,j}^l(m)\otimes
H_{i,m-j}^s(m+1)) \in\mathcal {T}_{2i+1,2i}\mathcal
{T}_{2i,2j+1}\subseteq \mathcal {T}_{2i+1,2j+1}.
$$
From (\ref{8}) we obtain
\begin{eqnarray}\label{13}
H_{i,j}^{j-l}(m)\otimes((s+1)H_{m-i,m-j}^{m-j-s-1}(m+1)+(i-s+1)H_{m-i,m-j}^{m-j-s}(m+1))\in\mathcal
{T}_{2i+1,2j+1}.
\end{eqnarray}
Since $l\in G_{m-i,j}(m)$ and $s\in G_{i,m-j}(m+1)$, from (\ref{13})
and Lemma \ref{24} we get
$$
H_{i,j}^{l'}(m)\otimes H_{m-i,m-j}^{s'}(m+1)\in\mathcal
{T}_{2i+1,2j+1},\quad l'\in G_{i,j}(m),  \; s'\in G_{m-i,m-j}(m+1),
$$
which implies $\mathcal {M}_{2i+1,2j+1}\subseteq\mathcal
{T}_{2i+1,2j+1}$.

\medskip
\textbf{Case 2.} $q=2j$ ($i+1\leq
j\leq\lfloor\frac{m}{2}\rfloor$).

By inductive hypothesis, we have
$$
H_{m-i,m-j}^l(m)\otimes H_{i,j}^s(m+1)\in\mathcal
{M}_{2i,2j}\subseteq\mathcal {T}_{2i,2j}, \quad l\in G_{m-i,m-j}(m),
\; s\in G_{i,j}(m+1).
$$
Thus by (\ref{12}) we obtain
$$
A_{2i+1,2i}(H_{m-i,m-j}^l(m)\otimes H_{i,j}^s(m+1))\in\mathcal
{T}_{2i+1,2i}\mathcal {T}_{2i,2j}\subseteq \mathcal {T}_{2i+1,2j}.
$$
From (\ref{8}), we have
\begin{eqnarray}\label{15}
H_{i,m-j}^{m-j-l}(m)\otimes((s+1)H_{m-i,j}^{j-s-1}(m+1)+(i-s+1)H_{m-i,j}^{j-s}(m+1))\in\mathcal
{T}_{2i+1,2j}.
\end{eqnarray}
Since $l\in G_{m-i,m-j}(m)$ and $s\in G_{i,j}(m+1)$, by (\ref{15})
and Lemma \ref{25} we have
$$
H_{i,m-j}^{l''}(m)\otimes H_{m-i,j}^{s''}(m+1)\in\mathcal
{T}_{2i+1,2j}, \quad l''\in G_{i,m-j}(m), \; s''\in G_{m-i,j}(m+1),
$$
which yields $\mathcal {M}_{2i+1,2j}\subseteq\mathcal
{T}_{2i+1,2j}$.

\medskip
\textbf{Step 2.2.}  Show $\mathcal {M}_{2i+2,q}\subseteq \mathcal
{T}_{2i+2,q}$ $(2i+2\leq q\leq m)$.

The proof of this step is similar to that of Step 2.1 and we omit it
here.

Hence the desired result follows by induction. $\qed$

\begin{thm}\label{17}
Let $\mathcal {T}$ be the Terwilliger algebra of the Odd graph
$O_{m+1}$ and $\mathcal {M}$ be the algebra defined in $(\ref{9})$.
Then $\mathcal {T}=\mathcal {M}$.
\end{thm}
\textbf{Proof.} Combining Lemma \ref{10} and Lemma \ref{20}, the
desired result follows. $\qed$

Since the generating matrices of each vector space in
(\ref{1})-(\ref{23}) are linearly independent, we have the following
result.

\begin{prop}\label{18}
The Terwilliger algebra $\mathcal {T}$ of the Odd graph $O_{m+1}$
has a basis:
$$
\{L(H_{m-i,m-j}^l(m)\otimes H_{i,j}^s(m+1)), \ l\in G_{m-i,m-j}(m),
\ s\in G_{i,j}(m+1)\}_{i,j=0}^m.
$$
\end{prop}

\begin{cor}\label{21}
The dimension of $\mathcal {T}$ is ${m+4 \choose 4}$.
\end{cor}
\textbf{Proof.} By Proposition \ref{18} we get
\begin{eqnarray*}
&&\dim\mathcal
{T}\\
&=&\sum_{i,j=0}^m|G_{m-i,m-j}(m)||G_{i,j}(m+1)|\\
&=&\sum_{i,j=0}^m(\min(m-i,m-j)-\max(0,m-i-j)+1)(\min(i,j)-\max(0,i+j-m-1)+1).
\end{eqnarray*}
By zigzag calculation, the desired result follows.
$\qed$

\section*{Acknowledgement} This research is partially supported by NSF of China, NCET-08-0052, and the Fundamental Research Funds for the Central
Universities of China.

\end{CJK*}

\end{document}